\documentclass[12pt]{amsart}
\newcommand{\bC}{{\mathbb C}}

\newcommand{\bP}{{\mathbb P}}

\newcommand{\bZ}{{\mathbb Z}}
\newcommand{\cA}{{\mathcal A}}

\newcommand{\cH}{{\mathcal H}}
\newcommand{\cI}{{\mathcal I}}
\newcommand{\cJ}{{\mathcal J}}
\newcommand{\cL}{{\mathcal L}}
\newcommand{\cO}{{\mathcal O}}

\newcommand{\cQ}{{\mathcal Q}}

\newcommand{\tC}{\widetilde{C}}

\newcommand{\tZ}{\widetilde{Z}}
\newcommand{\ra}{\rightarrow}
\newcommand{\lra}{\longrightarrow}

\newcommand{\inj}{\hookrightarrow}
\newcommand{\surj}{-\hspace{-4pt}\rightarrow\hspace{-19pt}{\rightarrow}\hspace{3pt}}

\newcommand{\D}{\Delta}
\newcommand{\T}{\Theta}

\newcommand{\Sg}{\mathfrak{S}}

\newtheorem{theorem}{Theorem}[section]
\newtheorem{lemma}[theorem]{Lemma}
\newtheorem{proposition}[theorem]{Proposition}
\newtheorem{corollary}[theorem]{Corollary}

\theoremstyle{definition}

\theoremstyle{remark}

\numberwithin{equation}{section}

\begin{document}

\title{Deforming curves in jacobians to non-jacobians I: \\
curves in $C^{ (2) }$}

\author{E. Izadi}
\address{Department of Mathematics, Boyd
Graduate Studies Research Center, University of Georgia, Athens, GA
30602-7403, USA}
\email{izadi@math.uga.edu}
\thanks{This material is based upon work partially supported by the
National Security Agency under Grant No. MDA904-98-1-0014 and the
National Science Foundation under Grant No. DMS-0071795. Any opinions,
findings and conclusions or recomendations expressed in this material
are those of the author and do not necessarily reflect the views of
the National Science Foundation (NSF) or the National Security Agency
(NSA).}





\maketitle

\section*{Introduction}

Jacobians of curves are the best understood abelian varieties. There
are many geometric ways of constructing curves in jacobians whereas it
is difficult to construct interesting curves in most other abelian
varieties. In this paper and its sequels we introduce methods for
determining whether a given curve in a jacobian deforms with it when
the jacobian deforms to a non-jacobian. We apply these methods to a
particular class of curves in a jacobian (see below). One of our
motivations is the Hodge-theoretic question of which multiples of the
minimal cohomology class on a given principally polarized abelian
variety can be represented by an algebraic curve (see \cite{I15} for
other motivations). In the known examples, the construction of such
curves often also leads to parametrizations of the theta divisor. Let
us begin by summarizing some of what is known about this question.

For a principally polarized abelian variety ({\em ppav}) $(A ,\Theta)$
of dimension $g$ over $\bC$, let $[\Theta]\in H^2 (A,\bZ)$ be the
cohomology class of $\Theta$. The class $\frac{[\Theta]^{ g-1 }}{ (g-1
)!}$ is called the minimal cohomology class for curves in $(A
,\Theta)$. We will assume $g\geq 4$, since otherwise all
multiples of the minimal class are represented by algebraic curves.

If $(A ,\Theta) = (JC := Pic^0 C,\Theta )$ is the jacobian of a
smooth, complete and irreducible curve $C$ of genus $g$, the choice of
any invertible sheaf $\cL$ of degree $1$ on $C$ gives an embedding of
$C$ in $JC$ via
\[
\begin{array}{rcl}
C &\inj & JC\\
p &\mapsto &\cO_C (p)\otimes\cL^{ -1}.
\end{array}
\]
Such a map is called an Abel map and the image of $C$ by it an Abel
curve. The cohomology class of an Abel curve is the minimal class
$\frac{\Theta^{ g-1 }}{ (g-1 )!}$. By a theorem of Matsusaka
\cite{matsusaka59}, the minimal class is represented by an algebraic
curve $C$ in $(A ,\Theta)$ if and only if $(A ,\Theta)$ is the
polarized jacobian $( JC ,\Theta )$ of $C$. The sum of $g-1$ Abel
curves in $JC$ is a theta divisor.

In Prym varieties even multiples of the minimal class are represented
by algebraic curves:

For an \'etale double cover $\pi :\tC\ra C$ of smooth curves ($C$ of
genus $g+1$), the involution $\sigma :\tC\ra\tC$ of the cover acts on
the jacobian $J\tC$ and the Prym variety $P$ of $\pi$ is defined as
\[
P := im(\sigma - 1)\subset J\tC.
\]
The intersection of a symmetric theta divisor of $J\tC$ with $P$ is
$2\Xi $ where the divisor $\Xi$ defines a principal polarization on
$P$. Via $\sigma -1$, the image of an Abel embedding of $\tC$ in $P$
gives an embedding of $\tC$ in $P$ and its image is called a
Prym-embedded curve. The class of a Prym-embedded curve is
$2\frac{[\Xi]^{ g-1 }}{ (g-1 )!}$. An interesting question is when, in
a Prym variety, are
any odd multiples of the minimal class represented by algebraic
curves.

A further generalization is the notion of Prym-Tjurin variety:

Suppose there is a generically injective map from a smooth complete
curve $X$ into a ppav $P$ with theta divisor $\Xi$ such that the class
of the image of $X$ is $m\frac{[\Xi]^{ g-1 }}{ (g-1 )!}$. This yields
a map $JX\ra P$ with transpose $P\ra JX$. We say $(P,\Xi)$ is a
Prym-Tjurin variety for $X$ if the map $P\ra JX$ is injective. This is
equivalent to the existence of an endomorphism $\phi :JX\ra JX$ such
that
\[
P = P(X,\phi) := im(\phi -1)\subset JX \qquad\hbox{and}\qquad
(\phi -1) (\phi -1 +m) = 0
\]
where the numbers denote the endomorphisms of $JX$ given by
multiplication by those numbers.

Welters showed that every principally polarized abelian variety is a
Prym-Tjurin variety. Birkenhacke and Lange showed that every
principally polarized abelian variety is a Prym-Tjurin variety for an
integer $m\leq 3^g (g-1)!$ (see \cite{birkenhakelange92} page 374
Corollary 2.4)\footnote{Their proof uses $3$-theta divisors. Using the
fact that a general $2$-theta divisor is smooth, the exact same proof
would give $m\leq 2^g (g-1)!$. For abelian varieties with a smooth
theta divisor, the same proof would give $m\leq (g-1)!$.}.

The set of integers $n$ such that $n\frac{[\Theta]^{ g-1 }}{ (g-1 )!}$
is represented by an algebraic curve (union $\{ 0\}$) is a
semi-subgroup of $\bZ$. There is therefore a unique minimal set of
positive integers $\{ m^A_0 < \ldots < m^A_{ r_A }\}$ such that
$n\frac{[\Theta]^{ g-1 }}{ (g-1 )!}$ is represented by an algebraic
curve if and only if either $n=m^A_i$ for some $i$ or $n > m^A_{ r_A
}$ is a multiple of $d_A$ where $d_A$ is the gcd of the $m^A_i$. These
various integers can be used to define stratifications of the moduli
space $\cA_g$ of ppav. The stratification associated to $m^A_0$ is
related to a stratification of $\cA_g$ by other geometric invariants
as discussed by Beauville \cite{beauville89} and Debarre
\cite{debarre92}. Debarre proved in \cite{debarre94} that $m^A_0$ is
at least $\sqrt{\frac{ g}{ 8}} -\frac{1}{4 }$ if $(A ,\Theta)$ is
general.

In this paper and the next we use first order obstructions to
deformations to identify certain families of curves which could
potentially deform to non-jacobians. Our method can be applied to
subvarieties $X$ of $JC$ contained in many translates of the theta
divisor. For an integer $e\geq 2$ let $C^{ (e)}$ be the $e$-th symmetric
power of $C$. Choose $a\in Pic^{ g-1 } C$ and define $\rho$ as the map
\[
\begin{array}{rccl}
\rho : & C^{ (g-1) } &\lra & JC\\
& D &\longmapsto &\cO_C (D)\otimes a^{ -1}
\end{array}
\]
whose image is a theta divisor, say $\T_a$. The idea is to use Green's
sequence \cite{green84} (see Section \ref{secttrans})
\[
0\lra T_{ C^{ (g-1 )}}\lra\rho^* T_{JC}\lra\cI_{ Z_{ g-1 }} (\Theta_a )\lra
0
\]
where $Z_{ g-1}$ is the locus where $\rho$ fails to be an embedding,
the letter $T$ denotes tangent sheaves and the letter $\cI$ ideal
sheaves. Given an infinitesimal deformation $\eta\in H^1 (T_{JC})$,
the curve $X$ deforms with $JC$ in the direction of $\eta$ if and only
if the image of $\eta$ by the first order obstruction map
\[
\nu : H^1 (T_{JC} )\lra H^1 (N_{ X/JC })
\]
where $N_{ X/JC }$ is the normal sheaf to $X$ in $JC$, is zero (see
Section \ref{sectdefobs}). As we
shall see in Section \ref{secttrans}, the map
\[
H^1 (T_{JC} )\lra H^1 (\cI_{ Z_{ g-1 }} (\Theta_a ) |_X)
\]
factors through $\nu$. It follows that if, for some $a$, the image of
$\eta$ in $H^1 (\cI_{ Z_{ g-1 }} (\Theta_a ) |_X)$ does not vanish,
then $\nu (\eta)$ does not vanish either. For the examples that we
have chosen (see below) we prove the stronger statement that the image
of $\eta$ by the map
\[
H^1 (T_{JC} )\lra H^1 (\cI_{ Z_{ g-1 }\cap X} (\Theta_a ))
\]
is not zero. This map factors into the composition (see Section
\ref{secttrans})
\[
H^1 (T_{JC} )\lra H^0 (\cO_{ Z_{ g-1 }\cap X} (\T))\lra H^1 (\cI_{ Z_{
g-1 }\cap X} (\Theta_a )).
\]
We analyze these two maps separately in Sections \ref{sectOmap} and
\ref{sectcobmap}. Section \ref{sectapp} (the Appendix) contains some
useful technical results.

The curves that we have chosen for the illustration of the above
method are the natural generalizations of smooth Prym-embedded curves
in tetragonal jacobians. More precisely, let $C$ be a curve of genus
$g$ with a $g^1_d$ (a pencil of degree $d$). We define curves $X_e
(g^1_d )$ whose reduced support is
\[
{X_e ( g^1_d )}_{ red } :=\{ D_e : \exists D\in C^{ (d-e) } \hbox{ such that }
D_e + D\in g^1_d\}\subset C^{ (e) }
\]
(for the precise scheme-theoretical definition see \ref{defX} when
$e=2$ and \cite{I15} for $e>2$). If $d\geq e+1$, then $X_e
(g^1_d)$ can be non-trivially mapped to $JC$ by subtracting a fixed
divisor of degree $e$ on $C$. Given a one-parameter infinitesimal
deformation of the jacobian of $C$ out of the jacobian locus $\cJ_g$
we ask when the curve $X_e ( g^1_d )$ deforms with it. In this paper
we prove the following
\begin{theorem}
Suppose $C$ non-hyperelliptic and $d\geq 4$. Suppose that the curve
$X := X_2 ( g^1_d )$ deforms in a direction $\eta$ out of $\cJ_g$. Then
\begin{enumerate}
\item either $d = 4$
\item or $d=5$, $h^0 (g^1_5) = 3$, $C$ has genus $4$, $L$ is
base-point-free, the curve $C$ has only one $g^1_3$ with a triple
ramification point $t$ such that $5t\in L\subset |K -t |$ and $X_2
(g^1_3)$ meets $X$ only at $2t$ with intersection multiplicity $4$.

\end{enumerate}
\end{theorem}

For $d=3$ the image of $X_2 (g^1_3 )$ in $JC$ is an Abel curve, hence
cannot deform out of $\cJ_g$ \cite{matsusaka59}. For $d=4$, the curve
$X_2 ( g^1_4 )$ is a Prym-embedded curve \cite{recillas74}, hence
deforms out of $\cJ_g$ into the locus of Prym varieties. For $d=5$,
$h^0 (g^1_d) = 3$ and $g=4$ (with only one $g^1_3$) or $g=5$ it is
likely that $X_2 ( g^1_5 )$ deforms out of $\cJ_g$. An interesting
question is to describe, as geometrically as possible, the
deformations of $( JC,\Theta )$ with which $X_2 ( g^1_5 )$ deforms.

For $e > 2$, the analogous result would be the following. The curve
$X_e (g^1_d )$ deforms out of $\cJ_g$ only if
\begin{itemize}
\item either $e= h^0 (g^1_d)$ and $d=2e$
\item or $e= h^0 (g^1_d ) -1$ and $d=2e+1$.
\end{itemize}
We prove this in \cite{I15} for $e\leq g-3$ under certain hypotheses
of genericity. This shortens the list of the families of curves whose
deformations we need to consider. For more details see \cite{I15}.

So we have some families of curves which could possibly deform to
non-jacobians. We need a different approach to prove that higher order
obstructions to deformations vanish: this will be presented in detail
in the forth-coming paper \cite{I16} and the idea behind it is the
following. For each $\T_a$ containing $X$, one has the map of
cohomology groups of normal sheaves
\[
H^1 (N_{ X/JC })\lra H^1 (N_{\T_a /JC } |_X ) = H^1 (\cO_X (\T_a ))
\]
whose kernel contains all the obstructions to the deformations of $X$
since we will only consider algebraizable deformations of $JC$ for
which the obstructions to deforming $\T_a$ vanish. If one can prove
that the intersection of these kernels is the image of the first order
algebraizable deformations of $JC$, i.e., the image of $S^2 H^1 (\cO_C
)\subset H^1 (T_{ JC})$, it will follow that the only obstructions to
deforming $X$ with $JC$ are the first order obstructions.

Finally, we would like to mention that from curves one can obtain
higher-dimensional subvarieties of an abelian variety. For a
discussion of this we refer the reader to \cite{I14}.

\section{Notation and Conventions}\label{notation}

We will denote linear equivalence of divisors by $\sim$.

For any divisor or coherent sheaf $D$ on a scheme $X$, denote by $h^i
(D)$ the dimension of the cohomology group $H^i (D) = H^i (X, D)$. For
any subscheme $Y$ of $X$, we will denote $\cI_{Y/X}$ the ideal sheaf
of $Y$ in $X$ and $N_{ Y/X }$ the normal sheaf of $Y$ in $X$. When
there is no ambiguity we drop the subscript ${}_X$ from $\cI_{Y/X}$ or
$N_{ Y/X }$. The tangent sheaf of $X$ will be denoted by $T_X :=\cH om
(\Omega^1_X ,\cO_X )$ and the dualizing sheaf of $X$ by $\omega_X$.

We let $C$ be a smooth non-hyperelliptic curve of genus $g$ over the
field $\mathbb C$ of complex numbers. For any positive integer $n$,
denote by $C^n$ the $n$-th Cartesian power of $C$ and by $C^{ (n) }$
the $n$-th symmetric power of $C$. We denote $\pi_n : C^n\ra C^{ (n)
}$ the natural quotient map. Note that $C^{ (n) }$ parametrizes the
effective divisors of degree $n$ on $C$.

We denote by $K$ an arbitrary canonical divisor on $C$. Since $C$ is
not hyperelliptic, its canonical map $C\ra |K|^*$ is an embedding and
throughout this paper we identify $C$ with its canonical image. For a
divisor $D$ on $C$, we denote by $\langle D\rangle$ its span in $| K
|^*$.

Since we will mostly work with the Picard group $Pic^{ g-1 }C$ of
invertible sheaves of degree $g-1$ on $C$, we put $A := Pic^{g-1}
C$. Let $\Theta$ denote the natural theta divisor of $A$, i.e.,
\[
\Theta :=\{\cL\in A : h^0 (\cL) > 0\}\; .
\]
The multiplicity of $\cL\in\Theta$ in $\Theta$ is $h^0 (\cL)$
(\cite{ACGH} Chapter VI p. 226). So the singular locus of $\Theta$ is
\[
Sing (\Theta ):=\{\cL\in A : h^0 (\cL )\geq 2\}\; .
\]
There is a map
\[
\begin{array}{rcl}
Sing_2 (\Theta ) & \lra & |I_2 (C)|\\
\cL & \longmapsto & Q (\cL ) :=\cup_{ D\in |\cL | }\langle D\rangle
\end{array}
\]
where $Sing_2 (\Theta )$ is the locus of points of order $2$ on
$\Theta$ and $| I_2 (C) |$ is the linear system of quadrics containing
the canonical curve $C$. This map is equal to the map sending $\cL$ to the
(quadric) tangent cone to $\Theta$ at $\cL$ and its image $\cQ$
generates $| I_2 (C) |$ (see \cite{green84} and
\cite{smithvarley90}). Any $Q (\cL )\in\cQ$ has rank $\leq 4$. The
singular locus of $Q (\cL)$ cuts $C$ in the sum of the base divisors
of $|\cL |$ and $ |\omega_C\otimes\cL^{ -1 } |$. The rulings of $Q$
cut the divisors of the moving parts of $|\cL |$ and $
|\omega_C\otimes\cL^{ -1 } |$ on $C$ (see \cite{andreottimayer67}).

For any divisor or invertible sheaf $a$ of degree $0$ and any
subscheme $Y$ of $A$, we let $Y_a$ denote the translate of $Y$ by
$a$. By a $g^r_d$ we mean a (not necessarily complete) linear system
of degree $d$ and dimension $r$. We call $W^r_d$ the subvariety of
$Pic^d C$ parametrizing invertible sheaves $\cL$ with $h^0 (\cL) > r$.

For any effective divisor $E$ of degree $e$ on $C$ and any positive
integer $n\geq e$, let $C^{( n-e) }_E\subset C^{ (n )}$ be the image
of $C^{ (n-e )}$ in $C^{ (n) }$ by the morphism $D\mapsto D+E$. Write
$C_t := C^{ (1)}_t$, and for any divisor $E =\sum_{i=1}^{ r} n_i t_i$ on
$C$, let $C_E$ denote the divisor $\sum_{ i=1 }^{ r} n_i C_{t_i}$ on
$C^{ (2) }$. For a linear system $L$ on $C$, we denote by $C_L$ any
divisor $C_E\subset C^{ (2) }$ with $E\in L$. We let $\delta$ denote
the divisor class on $C^{ (2) }$ such that $\pi_2^*\delta\sim\Delta$
where $\Delta$ is the diagonal of $C^2$.

By infinitesimal deformation we always mean {\em flat} first order
infinitesimal deformation.

\section{The curve $X$ and the first order obstruction map for its
deformations}\label{sectdefobs}

\subsection{}\label{defX} Let $L$ be a pencil of degree $d\geq 4$ on
$C$. Let $M$ be the moving part of $L$ and let $B$ be its base
divisor. Define the curve $X := X_2(L)$ as a divisor on $C^{ (2) }$ in
the following way
\[
X = X_2(L) := X_M + C_B\subset C^{ (2) }
\]
where $X_M := X_2(M)$ is the reduced curve
\[
X_M = X_2(M) :=\{ D_2 : \exists D\in C^{ (d-2) } \hbox{ such that }
D_2 + D\in M\}.
\]

\begin{lemma}
We have
\[
X\sim C_L -\delta
\]
and $X$ has arithmetic genus
\[
g_X =\frac{(d-2)(2g+d-3)}{2}.
\]
\end{lemma}
\begin{proof} Pull back to $C^2$, restrict to the fibers of the two
projections and use the See-Saw Theorem. For the arithmetic genus use
the exact sequence
\[
0\lra\cO_{ C^{ (2) }} (-X)\lra\cO_{ C^{ (2) }}\lra\cO_X\lra 0
\]
and the results of Appendix \ref{LELE}.
\end{proof}

\subsection{}\label{def}
Choose $g-3$ general points $p_1 ,\ldots , p_{g-3}$ in $C$ and embed
$C^{ (2) }$ in $C^{ (g-1) }$ and $A$ by the respective morphisms

\[
\begin{array}{llllllll}
C^{ (2) } &\lra & C^{ (g-1) } & & & C^{ (2) } &\lra & A\\
D_2 &\longmapsto & D_2 +\sum_{ i=1 }^{ g-3 } p_i &
& & D_2 &\longmapsto &\cO_C (D_2 +\sum_{ i=1 }^{ g-3 } p_i).
\end{array}
\]
Identify $X$ and $C^{ (2) }$ with their images by these maps.

Recall the usual exact sequence
\[
\cI_{ X/A }/\cI_{ X/A }^2\lra\Omega^1_A |_X\lra\Omega^1_X\lra 0.
\]
The curve $X$ is a local complete intersection scheme because it is a
divisor in $C^{ (2) }$. Using this, local calculations show that the
above sequence can be completed to the exact sequence
\[
0\lra\cI_{ X_{ red }/X }\cdot\cI_{ X/A }/\cI_{ X/A }^2\lra\cI_{ X/A
}/\cI_{ X/A }^2\lra\Omega^1_A |_X\lra\Omega^1_X\lra 0
\]
where $X_{ red }$ is the underlying reduced scheme of $X$. This
sequence can then be split into the following two short exact
sequences.
\begin{equation}\label{seqobs}
0\lra\cI_{ X_{ red }/X }\cdot\cI_{ X/A }/\cI_{ X/A }^2\lra\cI_{ X/A
}/\cI_{ X/A }^2\lra\cI_{ X/A }/\cI_{ X/A }^2 |_{ X_{ red }}\lra 0
\end{equation}
\[
0\lra\cI_{ X/A }/\cI_{ X/A }^2 |_{ X_{ red }}\lra\Omega^1_A
|_X\lra\Omega^1_X\lra 0
\]
from which we obtain the maps of exterior groups
\[
Ext^1 (\cI_{ X/A }/\cI_{ X/A }^2 |_{ X_{ red }} ,\cO_X )\lra H^1
(N_{ X/A })
\]
\[
H^1 ( T_A |_X )\surj Ext^1 (\cI_{ X/A }/\cI_{ X/A }^2 |_{ X_{ red
}} ,\cO_X ).
\]
The composition of the above two maps with restriction
\[
H^1 ( T_A )\lra H^1 ( T_A |_X )
\]
is the obstruction map
\[
\nu : H^1 (T_A )\lra H^1 ( N_{ X/A }).
\]
Given an infinitesimal deformation $\eta\in H^1 ( T_A )$, the curve
$X$ deforms with $A$ in the direction of $\eta$ if and only if $\nu
(\eta )= 0$ (see e.g. \cite{kollar96} Chapter 1 and
\cite{smithvarley90} for these deformation theory results). 

\subsection{} Using the fact that $X$ is a divisor in $C^{ (2) }$, a
local calculation shows that we have the usual exact sequence
\[
0\lra\cI_{ C^{ (2) }/A }/\cI_{ C^{ (2) }/A }^2 |_X\lra\cI_{ X/A }/\cI_{ X/A
}^2\lra\cI_{ X/C^{ (2) }}/\cI_{ X/C^{ (2) }}^2\lra 0
\]
whose dual is the exact sequence
\[
0\lra N_{ X/C^{ (2) }}\lra N_{ X/A }\lra N_{ C^{ (2) }/A } |_X\lra 0.
\]
From this we obtain the map
\[
H^1 ( N_{ X/A })\lra H^1 (N_{ C^{ (2) }/A } |_X)
\]
whose composition with $\nu$ we call $\nu_2$:
\[
\nu_2 : H^1 (T_A )\lra H^1 (N_{ C^{ (2) }/A } |_X).
\]
So, if $\nu (\eta )=0$, then, a fortiori, $\nu_2 (\eta )=0$.

\subsection{} The choice of the polarization $\Theta$ provides an
isomorphism $H^1 (T_A )\cong H^1 (\cO_C )^{\otimes 2}$ by which the
algebraic (i.e. globally unobstructed) infinitesimal deformations with
which $\Theta$ deforms are identified with the elements of the
subspace $S^2 H^1 (\cO_C )\subset H^1 (\cO_C )^{\otimes 2}\cong H^1
(T_A )$. Via this identification, the space of infinitesimal
deformations of $(A ,\Theta)$ as a jacobian is naturally identified
with $H^1 (T_C )\subset S^2 H^1 (\cO_C)$. The Serre dual of this last
map is multiplication of sections
\[
S^2 H^0 ( K)\lra H^0 ( 2K )
\]
whose kernel is the space $I_2 (C)$ of quadrics containing the
canonical image of $C$. Therefore, to say that we consider an
infinitesimal deformation of $(A ,\Theta )$ out of the jacobian locus,
means that we consider $\eta\in S^2 H^1 (\cO_C )\setminus H^1 (T_C )$
which is equivalent to say that we consider $\eta\in S^2 H^1 (\cO_C )$
such that there is a quadric $Q\in I_2 (C)$ with $(Q,\eta )\neq
0$. Here we denote by
\[
(,) : S^2 H^0 (K)\otimes S^2 H^1 (\cO_C )\lra S^2 H^1 (K)
\]
the pairing obtained from Serre Duality.

We fix such an infinitesimal deformation $\eta$ and prove that if $\nu_2
(\eta ) = 0$, then $d=4$ or $d=5$ and $h^0 (L )= 3$. For this
we use translates of $\Theta$ containing $C^{ (2) }$ and, a fortiori,
$X$.

\section{The translates of $\T$ containing $C^{ (2) }\supset X$ and
the first order obstruction map}\label{secttrans}

\begin{lemma} The surface $C^{ (2) }$ is contained in a translate
$\Theta_a$ of $\Theta$ if and only if there exists $\sum q_i\in C^{ (g-3) }$
such that $a =\sum p_i -\sum q_i$.
\end{lemma}
\begin{proof} For any points $q_1 ,\ldots ,q_{g-3}$ of $C$, the image of
$C^{ (2) }$ in $A$ is contained in the divisor $\Theta_{\sum p_i
  -\sum q_i }$. Conversely, if $C^{ (2) }$ is contained in a translate
$\Theta_a$ of $\Theta$, then we have $h^0 ( D_2 +\sum p_i -a) > 0$,
for all $D_2\in C^{ (2) }$. Equivalently, for all $D_2\in C^{ (2) }$,
we have $h^0 (K + a -\sum p_i -D_2 ) > 0$, i.e., $h^0 (K + a -\sum p_i
)\geq 3$ and $-a+\sum p_i$ is effective.
\end{proof}

\subsection{}
Choose $\cO_C (a)\in Pic^0 C$ such that $C^{ (2) }\subset\Theta_{-a}$
(i.e., $-a =\sum p_i -\sum q_i$ as above). Then $C^{ (2)
}_a\subset\Theta$. Let $\rho : C^{ (g-1 ) }\ra\Theta$ be the
natural morphism. Then (see \cite{green84} (1.20) p. 89) we have the
exact sequence
\begin{equation}\label{TCg-1}
0\lra T_{ C^{ (g-1 )}}\lra\rho^* T_A\lra\cI_{ Z_{ g-1 }} (\Theta )\lra
0
\end{equation}
where the leftmost map is the differential of $\rho$ and $Z_{ g-1 }$ is the subscheme of $C^{ (g-1 )}$ where $\rho$
fails to be an isomorphism. For the convenience of the reader we
mention that the scheme $Z_{ g-1 }$ is a determinantal scheme of
codimension $2$. If $g\geq 5$ or if $g=4$ and $C$ has two distinct
$g^1_3$'s, the scheme $Z_{ g-1}$ is reduced and is the
scheme-theoretical inverse image of the singular locus of
$\Theta$.

Combining sequence (\ref{TCg-1}) with the tangent bundles sequences for
$C^{ (2) }_a$, we obtain the commutative diagram with exact rows and
columns
\[
\begin{array}{ccccccccc}
& & 0 & & 0 & & & & \\
& & \downarrow & & \downarrow & & & & \\
& & T_{ C^{ (2) }_a } & = & T_{ C^{ (2) }_a } & & & & \\
& & \downarrow & &\downarrow & & & & \\
 &  & T_{ C^{ (g-1) }} |_{ C^{ (2) }_a } &\ra & T_A |_{ C^{ (2) }_a }
&\ra & \cI_{ Z_{ g-1 }} (\Theta ) |_{ C^{ (2) }_a } &\ra & 0 \\
& & \downarrow & &\downarrow & & || & & \\
 &  & N_{  C^{ (2) }_a/C^{ (g-1) }} &\ra &
N_{  C^{ (2) }_a/A } &\ra & \cI_{ Z_{ g-1
}} (\Theta ) |_{ C^{ (2) }_a } &\ra & 0 \\
& & \downarrow & &\downarrow & & & & \\
& & 0 & & 0 & & & & \\
\end{array}
\]
where the leftmost horizontal maps are injective if and only if $h^0
(\sum q_i )=1$. The restriction of this to $X_a$ gives the
commutative diagram with exact rows and columns
\[
\begin{array}{ccccccc}
T_{ C^{ (2) }_a } |_{ X_a } & = & T_{ C^{ (2) }_a } |_{ X_a } & & & & \\
\downarrow & &\downarrow & & & & \\
T_{ C^{ (g-1) }} |_{ X_a } &\ra & T_A |_{ X_a }
&\ra & \cI_{ Z_{ g-1 }} (\Theta ) |_{ X_a } &\ra & 0 \\
\downarrow & &\downarrow & & || & & \\
N_{  C^{ (2) }_a/C^{ (g-1) }} |_{ X_a } &\ra &
N_{  C^{ (2) }_a/A } |_{ X_a } &\ra & \cI_{ Z_{ g-1
}} (\Theta ) |_{ X_a } &\ra & 0 \\
\downarrow & &\downarrow & & & & \\
0 & & 0 & & & & \\
\end{array}
\]
whose cohomology gives the commutative diagram with exact rows and columns
\[
\begin{array}{ccccccc}
H^1 ( T_{ C^{ (2) }_a } |_{ X_a } ) & = & H^1 ( T_{ C^{ (2) }_a } |_{
X_a } ) & & & & \\
\downarrow & &\downarrow & & & & \\
H^1 ( T_{ C^{ (g-1) }} |_{ X_a } ) &\ra & H^1 ( T_A |_{ X_a } )
&\ra & H^1 ( \cI_{ Z_{ g-1 }} (\Theta ) |_{ X_a } ) &\ra & 0 \\
\downarrow & &\downarrow & & || & & \\
H^1 ( N_{  C^{ (2) }_a/C^{ (g-1) }} |_{ X_a } ) &\ra &
H^1 ( N_{  C^{ (2) }_a/A } |_{ X_a } ) &\ra & H^1 ( \cI_{ Z_{ g-1
}} (\Theta ) |_{ X_a } ) &\ra & 0 \\
\downarrow & &\downarrow & & & & \\
0 & & 0 & & & & .\\
\end{array}
\]
Therefore we have the commutative diagram
\[
\begin{array}{ccccc}\label{diagI}
H^1 ( T_A ) &\lra & H^1 (T_A |_{X_a} ) &\lra & H^1 (\cI_{ Z_{ g-1 }}
(\Theta ) |_{X_a} ) \\
\bigcup & \searrow & \downarrow & & || \\
S^2 H^1 (\cO_C ) & \lra & H^1 ( N_{ C^{ (2) }_a/A } |_{X_a} ) & \lra &
H^1 (\cI_{ Z_{ g-1 }} (\Theta ) |_{X_a} ).
\end{array}
\]
Translation by $a$ induces the identity on $H^1 (T_A )$ and isomorphisms
\[
H^1 (T_A |_{X_a} )\cong H^1 (T_A |_{X} )\qquad H^1 ( N_{ C^{ (2) }_a/A
} |_{X_a} )\cong H^1 ( N_{ C^{ (2) }/A } |_{X} )
\]
so that the kernel of
\[
\nu_2 : S^2 H^1 (\cO_C )\lra H^1 (N_{ C^{ (2) }/A } |_{X} )
\]
is equal to the kernel of the map
\[
S^2 H^1 (\cO_C )\lra H^1 (N_{ C^{ (2) }_a/A } |_{X_a} )
\]
obtained from $\nu_2$ by translation. Therefore the previous diagram proves
the following theorem.

\begin{theorem}\label{kernu}
The kernel of the map
\[
\nu_2 : S^2 H^1 (\cO_C )\lra H^1 (N_{ C^{ (2) }/A } |_{X} )
\]
is contained in the kernel of the map obtained from the above
\[
S^2 H^1 (\cO_C )\lra H^1 (\cI_{ Z_{ g-1 }} (\Theta ) |_{X_a} )
\]
for all $a$ such that $\Theta_{-a}$ contains $C^{ (2) }$.
\end{theorem}

\subsection{} We shall prove that for any $\eta\in S^2 H^1 (\cO_C
)\setminus H^1 (T_C)$, there exists $a$ such that $\Theta_{ -a }$
contains $C^{ (2) }$ and the image of $\eta$ by the map
\[
S^2 H^1 (\cO_C )\lra H^1 (\cI_{ Z_{ g-1 }} (\Theta ) |_{X_a} )
\]
is nonzero unless either $d=4$ or $d=5$, $h^0 (L) =3$ and $C$ has
genus $5$ or genus $4$ and only one $g^1_3$.

\subsection{}\label{prepquad} The latter map is the composition of
\[
S^2 H^1 (\cO_C )\lra H^1 (\cI_{ Z_{ g-1 }} (\Theta ))
\]
with restriction
\[
H^1 (\cI_{ Z_{ g-1 }} (\Theta ))\lra H^1 (\cI_{ Z_{ g-1 }}
(\Theta ) |_{X_a} ).
\]
From the natural map
\[
\cI_{ Z_{ g-1 }} (\T) |_{ X_a }\lra\cI_{ Z_{ g-1 }\cap X_a } (\T)
\]
we obtain the map
\[
H^1 (\cI_{ Z_{ g-1 }} (\Theta ) |_{X_a})\lra H^1 (\cI_{ Z_{ g-1 }\cap
{X_a}} (\Theta )).
\]
Therefore we look at the kernel of the composition
\begin{equation}\label{comp1}
\begin{split}
S^2 H^1 (\cO_C )\lra H^1 (\cI_{ Z_{ g-1 }} (\Theta ))
\lra H^1 (\cI_{ Z_{ g-1 }} (\Theta ) |_{ X_a })\\
\lra H^1 (\cI_{ Z_{ g-1 }\cap X_a} (\Theta )).
\end{split}
\end{equation}
From the usual exact sequence
\[
0 \lra \cI_{ Z_{ g-1 }} (\Theta )\lra \cO_{ C^{ g-1 }} (\Theta
)\lra \cO_{ Z_{ g-1 }} (\Theta )\lra  0,
\]
we obtain the embedding
\[
H^0 (\cO_{ Z_{ g-1 }} (\Theta ))\inj H^1 (\cI_{ Z_{ g-1 }} (\Theta )).
\]
By \cite{green84} p. 95, the image of $S^2 H^1 (\cO_C )$ in $H^1
(\cI_{ Z_{ g-1 }} (\Theta ))$ is contained in $H^0 (\cO_{ Z_{ g-1 }}
(\Theta ))$.
Now, using the commutative diagram with exact rows
\[
\begin{array}{ccccccccc}
0 &\lra &\cI_{ Z_{ g-1 }} (\Theta )&\lra &\cO_{ C^{ g-1 }} (\Theta
)&\lra &\cO_{ Z_{ g-1 }} (\Theta )&\lra & 0\\
& &\downarrow & &\downarrow & &\downarrow & & \\
0 &\lra &\cI_{ Z_{ g-1 }\cap X_a} (\Theta )&\lra &\cO_{X_a} (\Theta
)&\lra &\cO_{ Z_{ g-1 }\cap X_a} (\Theta )&\lra & 0,\\
\end{array}
\]
Composition \ref{comp1} is equal to the composition
\[\begin{split}
S^2 H^1 (\cO_C )\lra H^0 (\cO_{ Z_{ g-1 }} (\Theta ))\lra 
H^0 (\cO_{ Z_{ g-1 }\cap {X_a}} (\Theta ))\lra\\
\stackrel{ coboundary }{\lra}  H^1 (\cI_{ Z_{ g-1 }\cap {X_a}} (\Theta )).\end{split}
\]
By \cite{green84} p. 95 again, the first map is the following
\[
\begin{array}{ccc}
S^2 H^1 (\cO_C ) &\lra & H^0 (\cO_{ Z_{ g-1 }} (\Theta ))\\
\sum a_{ij}\frac{\partial^2 }{\partial z_i\partial z_j} &\longmapsto
&\sum a_{ij}\frac{\partial^2\sigma }{\partial z_i\partial z_j} |_
{ Z_{ g-1 }},
\end{array}
\]
where $\{ z_i\}$ is a system of coordinates on $A$ and $\sigma$ is a
theta function with divisor of zeros equal to $\Theta$.
So we have
\[
\begin{array}{cccccc}
S^2 H^1 (\cO_C ) &\lra & H^0 (\cO_{ Z_{ g-1 }} (\Theta )) &\lra &
H^0 (\cO_{ Z_{ g-1 }\cap X_a} (\Theta )) &\lra\\
\sum a_{ij}\frac{\partial^2 }{\partial
z_i\partial z_j} &\longmapsto &\sum a_{ij}\frac{\partial^2\sigma
}{\partial z_i\partial z_j} |_{ Z_{ g-1 }} &\longmapsto &\sum
a_{ij}\frac{\partial^2\sigma }{\partial z_i\partial z_j} |_{ Z_{
g-1 }\cap X_a} &\longmapsto \\
& & &\stackrel{coboundary}{\lra}
& H^1 (\cI_{ Z_{ g-1 }\cap X_a} (\Theta )) &\\
& & &\longmapsto & ? & .
\end{array}
\]
We will investigate the kernel of the composition of the first two
maps and that of the coboundary map separately. The kernel of the
composition of the first two maps is contained in (with equality if
and only if $Z_{ g-1 }\cap X_a$ is reduced) the annihilator of the
quadrics of rank $\leq 4$ which are the tangent cones to $\Theta$ at
the points of $\rho (Z_{ g-1 })\cap X_a = Sing (\Theta )\cap X_a$.

\section{The kernel of the map $S^2 H^1 (\cO_C )\ra H^0 (\cO_{ Z_{ g-1
}\cap X_a} (\Theta ))$}\label{sectOmap}

\subsection{}\label{sectdefZ} Let
$\tZ (X)\subset C^{ (g-3) }\times X$ be the closure of the subvariety
parametrizing pairs $(\sum q_i , D_2 )$ such that $h^0 (\sum q_i ) =
1$ (this is the case generically on $C^{ (g-3) }$ because dim$W^1_{
g-3 }\leq g-3 -2-1 = g-6$ by \cite{mumford74} pp. 348-350) and $h^0 (
D_2 +\sum q_i ) = 2$. Let $Z(X)\subset C^{ (g-3) }$ be the image of
$\tZ (X)$ by the first projection. Denote by $\tZ (X_M)$ and $Z (X_M)$ the 
corresponding objects for $X_M$. It follows from Corollary
\ref{tZonto} below that $\tZ (X)$ and $Z (X)$ are not empty and $\tZ
(X_M)$ and $Z (X_M)$ are not empty when either  the degree of $M$ is at least
$4$ or $g\geq 5$.

Given an infinitesimal deformation $\eta\in S^2 H^1 (\cO_C )\setminus
H^1 (T_C)$ we shall prove that there is always a component of $\tZ
(X)$ such that for $(\sum q_i, D_2 )$ general in that component the
tangent cone to $\Theta$ at $\cO_C (D_2 +\sum q_i )$ does not vanish
on $\eta$. This will follows from Corollary \ref{tZonto} below, given
that the image $\cQ$ of $Sing_2 (\Theta )$ in $|I_2 (C) |$ generates $|
I_2 (C) |$. By our remarks above, it implies a fortiori that the image
of $\eta$ in $H^0 (\cO_{ Z_{ g-1 }\cap X_a} (\Theta ))$ is nonzero for
$-a =\sum p_i -\sum q_i$. We will then show that generically on any
component of $\tZ (X)$ the coboundary map
\[
H^0 (\cO_{ Z_{ g-1 }\cap X_a} (\Theta )) \stackrel{coboundary}{\lra}
H^1 (\cI_{ Z_{ g-1 }\cap X_a} (\Theta ))
\]
is injective unless either $d=4$ or $d=5$ and $h^0 (L) =3$. It will
follow by Theorem \ref{kernu} that $\nu_2 (\eta )\neq 0$, hence $\nu
(\eta )\neq 0$ and $X$ does not deform with $\eta$ unless either $d=4$
or $d=5$ and $h^0 (L) =3$. We begin by computing the dimensions of
$Z(X)$ and $\tZ (X)$ and showing that $\tZ (X)$ maps onto $Sing
(\Theta)$.

\begin{lemma}
For any quadric $Q$ containing $C$, there exists $D_2\in X$ such that
$\langle D_2\rangle\subset Q$. If either $g\geq 5$, or the degree of $M$ is
at least $4$, then such a $D_2$ can be chosen to be in $X_M\subset X$.
\end{lemma}
\begin{proof} By Appendix \ref{appI2} the space $I_2 (C)$ can be
identified with $H^0 ( C^{ (2) }, C_K - 2\delta )$. So $Q$ corresponds
to a section $s_Q\in H^0 ( C^{ (2) }, C_K - 2\delta )$. The support of
the divisor $E_Q$ of zeros of $s_Q$ is the set of divisors $D\in C^{
(2) }$ such that $\langle D\rangle\subset Q$. So our lemma is
equivalent to saying that the intersection $E_Q\cap X$ is not
empty. This follows from the following computation of the intersection
number of $E_Q$ and $X$ where we use $d\geq 4$ and $g\geq 4$.
\[\begin{split}
X\cdot E_Q = (C_L -\delta )\cdot (C_K - 2\delta ) = d (2g-2) - 2d -
(2g-2) - 2 (g-1)\\
= d (2g-4) - 4g + 4\geq 4 (2g-4) - 4g + 4 = 4g - 12\geq 4.\end{split}
\]
The analogous calculation with $X_M$ instead of $X$ proves the second
asertion.
\end{proof}

\begin{lemma}
For any $g^1_{ g-1 }$ on $C$, there exists $D_2\in X$ such that $h^0 (
g^1_{ g-1 } - D_2 )> 0$. If either $g\geq 5$, or the degree of $M$ is
at least $4$, then such a $D_2$ can be chosen to be in $X_M\subset X$.
\end{lemma}
\begin{proof} This follows from the positivity of
the intersection number of $X$ and $X_2( g^1_{ g-1 })$:
\[\begin{split}
X\cdot X_2( g^1_{ g-1 }) = (C_L -\delta )\cdot (C_{ g^1_{ g-1 }} -\delta ) =
d (g-1) - d - (g-1) - (g-1)\\
= d (g-2) - 2g + 2\geq 4 (g-2) - 2g + 2 = 2g - 6\geq 2.\end{split}
\]
Tha analogus calculation with $X_M$ instead of $X$ proves the
second assertion.
\end{proof}

\begin{corollary}\label{tZonto}
The variety $\tZ(X)$ maps onto $Sing (\Theta )$ by $\rho$. If either
$g\geq 5$, or the degree of $M$ is at least $4$, then $\tZ(X_M)$ also
maps onto $Sing (\Theta )$ by $\rho$.
\end{corollary}
\begin{proof} It is sufficient to prove that the map $\tZ(X)\ra Sing
(\Theta)$ is dominant. A general point of $Sing (\Theta )$ is a
complete $g^1_{ g-1 }$ on $C$. By the previous lemma, there is a
divisor of $g^1_{ g-1 }$ which contains $D_2$. So there is $\sum
q_i\in C^{ (g-3) }$ with $D_2 +\sum q_i\in g^1_{ g-1 }$. The pair
$(\sum q_i , D_2)\in C^{ (g-3) }\times X$ maps to $g^1_{ g-1 }$ by
$\rho$. To see that $(\sum q_i , D_2)\in\tZ (X)$ for a general choice
of $g^1_{ g-1}$, it is sufficient to prove that $h^0 (\sum q_i ) =1$
for a general choice of $g^1_{ g-1}$.

If $g^1_{ g-1}$ is base-point-free, then this is automatic. If $g^1_{
g-1}$ has base points, then it is sufficient to prove that no divisor
$D_2\in X$ is contained in its base divisor. It follows from a theorem
of Mumford (see \cite{ACGH} p. 193) that a general $g^1_{ g-1 }$ is
base-point-free unless $C$ is either trigonal, bielliptic or a smooth
plane quintic. In all these cases, the base divisor can be chosen to
be general so that it contains no divisor $D_2\in X$.

The assertion about $\tZ( X_M )$ is proved similarly, using the
corresponding statements for $X_M$.
\end{proof}

\begin{proposition}\label{dimZ}
The varieties $\tZ(X)$ and $\tZ (X_M )$ (when non-empty) are
everywhere of dimension $\geq g-4$.
\end{proposition}
\begin{proof} By Corollary \ref{tZonto}, the variety $\tZ (X)$ is not
empty. To see that the dimensions of $\tZ (X)$ and $\tZ (X_M )$ are
everywhere $\geq g-4$, note that $h^0 (D_2 +\sum q_i)\geq 2$ is
equivalent to $D_2 +\sum q_i\in Z_{ g-1 }$. Requiring $D_2\in X$
(resp. $X_M$) imposes at most one condition on the pair $(\sum q_i
,D_2 )$. Since the dimension of $Z_{ g-1 }$ is $g-3$ \cite{green84},
the proposition follows.

\end{proof}

\subsection{}\label{Qeta}
Since quadrics associated to $g^1_{ g-1 }$'s generate $I_2 (C)$ (see
\cite{green84} and \cite{smithvarley90}), for any direction $\eta\in
S^2 H^1 (\cO_C )\setminus H^1 (T_C)$ there exists an irreducible
component $\cQ (\eta )$ of $\cQ$ such that for $Q$ general in $\cQ
(\eta )$, the quadric $Q$ is nonzero on $\eta$ (in fact $\cQ$ is
almost always irreducible but we do not need to go into this). Let
$\tZ (\eta )$ be an irreducible component of $\tZ (X)$ which maps onto
$\cQ (\eta )$ and let $Z (\eta )$ be the image of $\tZ (\eta )$ in
$C^{ (g-3) }$. If the degree of $M$ is at least $4$ or if $g\geq 5$,
choose $\tZ (\eta )$ and $Z (\eta )$ to be in $\tZ (X_M)$ and $Z
(X_M)$ respectively. Then, for $\sum q_i$ general in $Z (\eta )$, the image
of $\eta$ in the corresponding $H^0 (\cO_{ Z_{ g-1 }\cap X_a} (\Theta
))$ is nonzero.

\section{The coboundary map $H^0 (\cO_{ Z_{ g-1 }\cap X_a} (\Theta ))
\lra H^1 (\cI_{ Z_{ g-1 }\cap X_a} (\Theta ))$}\label{sectcobmap}

\begin{lemma}\label{noninj}
Suppose $\sum q_i\in Z(X)$ satisfies $h^0 (\sum q_i ) =1$. If the
coboundary map
\[
H^0 (\cO_{ Z_{ g-1 }\cap X_a} (\Theta )) \lra H^1 (\cI_{ Z_{ g-1
}\cap X_a} (\Theta ))
\]
is not injective, then
\[
H^0 (C , K -\sum q_i - L)\neq 0.
\]
\end{lemma}
\begin{proof}
Using the exact sequence
\[\begin{split}
0\lra H^0 (\cI_{ Z_{ g-1 }\cap X_a} (\Theta )) \lra H^0 (\cO_{X_a}
(\Theta )) \lra
\\
\lra H^0 (\cO_{ Z_{ g-1 }\cap X_a} (\Theta )) \lra H^1
(\cI_{ Z_{ g-1 }\cap X_a} (\Theta )),\end{split}
\]
we need to understand the sections of $\cO_{X_a} (\Theta )$ which
vanish on $Z_{ g-1 }\cap X_a$. Equivalently, translating everything by
$-a$, we need to understand the sections of $\cO_X (\Theta_{-a} )$
which vanish on $(Z_{ g-1 })_{-a}\cap X$. For this, we use the
embedding of $X$ in $C^{(2)}$:
\[
0\lra\cO_{ C^{ (2) }} (\Theta_{-a} - X)\lra\cO_{ C^{ (2) }}
(\Theta_{-a} )\lra\cO_X (\Theta_{-a} )\lra 0.
\]
Since $\cO_{ C^{ (2) }} (\Theta_{-a} -X )\cong\cO_{ C^{ (2) }} (C_{ K
-\sum q_i -L}) =\cL_{ 2, K -\sum q_i -L}$ and $\cO_{ C^{ (2) }}
(\Theta_{-a} )\cong\cO_{ C^{ (2) }} (C_{ K -\sum q_i } -\delta
)=\cL'_{ 2, K -\sum q_i }$ (see \ref{LELE} for this notation), by
Appendix \ref{LELE} this gives the exact sequence of cohomology
\[\begin{split}
0\lra S^2 H^0 (C, K-\sum_{ i=1 }^{ g-3 } q_i -L)\lra\wedge^2 H^0 (C, K
-\sum_{ i=1 }^{ g-3 } q_i)
\\
\lra H^0 (X,\Theta_{-a} )
\lra H^0 (C, K -\sum_{ i=1 }^{ g-3 } q_i -L )\otimes H^1 (C, K -\sum_{
i=1 }^{ g-3 } q_i -L ).
\end{split}
\]
Since $h^0 (\sum q_i ) = 1$, by Appendix \ref{ontoC2} the elements of
$H^0 ( C^{ (2) },\Theta_{-a} ) =\wedge^2 H^0 (C, K -\sum_{ i=1 }^{ g-3
} q_i)$ all vanish on $(Z_{ g-1 })_{ -a}\cap C^{ (2) }$, hence they
also vanish on $(Z_{ g-1})_{ -a}\cap X$. So if the coboundary map is
not injective, then there must be elements of $H^0 (X,\Theta_{-a} )$
which are not restrictions of elements of $H^0 ( C^{ (2) },\Theta_{-a}
)$. In particular, by the above exact sequence, we must have $H^0 (C,
K-\sum_{ i=1 }^{ g-3 } q_i -L)\neq 0$.
\end{proof}

For $\eta\in S^2 H^1 (\cO_C )\setminus H^1 (T_C )$, define $\cQ (\eta
)$, $\tZ (\eta )$, $Z (\eta )$ as in Paragraph \ref{Qeta}. We have

\begin{theorem}
Suppose $X$ deforms infinitesimally with $A$ in a direction $\eta\in
S^2 H^1 (\cO_C )\setminus H^1 (T_C )$, then either $d=4$ or $d=5$ and
$h^0 (L )=3$. Furthermore, in this case we can choose $Z(\eta )$ to be
of dimension at least $g-4$ unless $g\geq 7$, $C$ is a double cover of
a smooth curve of genus $2$ and $L$ is the inverse image of the
$g^1_2$ on the curve of genus $2$.
\end{theorem}
\begin{proof} Let $\sum q_i\in Z(\eta )$ be general so that in
particular we have $h^0 (\sum q_i )=1$ (see \ref{sectdefZ}). Then, as
we noted in \ref{Qeta}, the image $\overline{\eta}$ of $\eta$ in $H^0
(\cO_{ Z_{ g-1 }\cap X_a} (\Theta ))$ is not zero. Since $X$ deforms
with $\eta$, by Theorem \ref{kernu}, the image of $\eta$ in $H^1
(\cI_{ Z_{ g-1 }} (\Theta ) |_{X_a} )$ is zero. So
$\overline{\eta}$ is in the kernel of the coboundary map
\[
H^0 (\cO_{ Z_{ g-1 }\cap X_a} (\Theta ))\lra H^1 (\cI_{ Z_{ g-1 }}
(\Theta ) |_{X_a} )
\]
which is therefore not injective. It follows, by Lemma \ref{noninj},
that $h^0 (K -\sum q_i -L ) >0$.

Since the dimension of $\tZ (\eta )$ is at least $g-4$ (see
Proposition \ref{dimZ}) and $X$ is one-dimensional, the dimension of
$Z (\eta )$ is at least $g-5$. If the genus is $4$, then since $\tZ
(\eta )$ is not empty, the dimension of $Z (\eta)$ is at least $g-4$.

If $Z(\eta)$ has dimension $\geq g-4$, then by the above discussion we
have $h^0 (K -\sum q_i -L ) >0$ for a $(g-4)$-dimensional family of
$\sum_{ i=1 }^{ g-3 } q_i$ (in $Z(\eta)$). So $h^0 (K-L)\geq g-3$ and,
by Clifford's Lemma, since $C$ is not hyperelliptic, we have $2(
g-3-1) < 2g-2-d$ or $d\leq 5$. If $d=5$, then clearly $h^0 (L) =3$.

Suppose now that every component $Z(\eta)$ has dimension $g-5$. Then
$g\geq 5$ by the above and $Z(\eta )\subset Z (X_M)$. Here Clifford's
Lemma only gives us $d\leq 7$ so we use the following argument. Since
$h^0 (\sum q_i ) =1$ generically on $Z (\eta )$, the $|\sum q_i |$
form a $(g-5)$-dimensional family of linear systems and so do the $|K
-\sum q_i |$. Writing $|K -\sum q_i | = L + B'$, the $B'$ vary in a
family of effective divisors of dimension $\geq g-5$. Therefore the
degree of $B'$ is at least $g-5$ and $d + g-5\leq 2g-2 -(g-3)$, i.e.,
$d\leq 6$. Next $\tZ (\eta)$ has dimension $g-4$ and the general
fibers of $\tZ (\eta )\ra Z(\eta )$ are one-dimensional, all equal to
a union of components of $X_M$, say $X'$. Since we can suppose $h^0
(\sum q_i) =1$ (see \ref{sectdefZ}), the condition $h^0 (D_2 +\sum q_i
)\geq 2$ for all $D_2\in X'$ means that $\langle\sum
q_i\rangle\cap\langle D_2\rangle\neq\emptyset$ for all $D_2\in
X'$. Therefore the projection of center $\langle\sum q_i\rangle$ from
the canonical curve $C$ to $\bP^2 = |K -\sum q_i |^*$ is not
birational to its image. So there is a nonconstant map $\kappa :C\ra
C'$ of degree $\geq 2$ with $C'$ smooth such that
\[
X'\subset\{ D_2\in C^{ (2) } : \exists t\in C' \hbox{ such that }
D_2\leq\kappa^* t\}.
\]
and
\[
|K-\sum q_i| =\kappa^* N + B_0
\]
where $N$ is a two-dimensional linear
system $N$ on $C'$ and $B_0$ is the base divisor of $|K -\sum q_i |$.

Consider now the linear systems $| K -\sum q_i |$. As $\sum q_i$
varies in $Z (\eta)$, the divisors of these linear systems form a
$(g-3)$-dimensional family of divisors. Therefore we have
\[
g-3\leq deg (B_0 ) + deg (N)\qquad\hbox{or}\qquad deg (B_0)\geq g-3-deg(N).
\]
Combining this with the equality
\[
deg (B_0 ) + deg (\kappa )deg (N) = g+1
\]
we obtain
\[
deg (N) (deg (\kappa ) -1 )\leq 4.
\]
Since $N$ has degree at least $2$ we first obtain
\[
deg (\kappa )\leq 3.
\]

If $deg (\kappa ) = 3$, then $deg (N) =2$ and $C'$ is a conic in
$\bP^2$. Hence $C$ is trigonal and $X' = X_2 (g^1_3 )=\{ D_2 : h^0 (
g^1_3 - D_2 )>0\}$. In this case $Z(X) = C^{ (g-3) }$ since for any
$\sum q_i\in C^{ (g-3) }$, if we take $D_2 = g^1_3 - q_1\in X'$, then
$h^0 ( D_2 -\sum q_i )\geq 2$. So $Z (\eta ) = Z(X)$ is of dimension
$g-3$ which is contrary to our hypothesis.

Therefore $\kappa$ has degree $2$ and $X' =\{\kappa^* t: t\in
C'\}$. In this case, since $C$ is not hyperelliptic, $C'$ is
birational to a plane curve of degree $3$ or $4$ and has genus $1,2$
or $3$. If $C'$ is elliptic, then any divisor $\sum_{ i=1 }^{ g-5 }
q_i +\kappa^* q$ is in $Z (\eta)$ and $Z (\eta )$ is of dimension
$g-4$ which is against our hypothesis. If $C'$ has genus $\geq 2$,
then its plane model has degree $4$. If $C'$ has genus $3$, then it
has only one $g^2_4$ which is then $N$. This implies that $|
K-\kappa^* N|$ has dimension $\geq g-5$ since $h^0 (K-\kappa^* N -\sum
q_i ) > 0$ for $\sum q_i$ in a $(g-5)$-dimensional family of effective
divisors. Therefore $h^0 (\kappa^* N )\geq 5$ by the Riemann-Roch
Theorem. However, this is impossible as $|K-\sum q_i| = B_0 +\kappa^*
N$ is a complete linear system of dimension $2$ for a general $\sum
q_i\in Z(\eta)$.

So $C'$ has genus $2$ and its plane model has a double point: $N=
g^1_2 + t_1 +t_2$ for some points $t_1$ and $t_2$ on $C'$ such that
$t_1 + t_2\not\in g^1_2$. We obtain $deg(B_0) = g-7$ and, for $\sum
q_i\in Z(\eta )$ general, $B_0$ is a general effective divisor of
degree $g-7$ on $C$. In particular, $g\geq 7$. Furthermore, since
$B_0$ is general and $h^0 (\kappa^* N + B_0 -L ) >0$ for all $B_0$, we
obtain
\[
h^0 (\kappa^* N -L )>0.
\]
Now, since the $|K -\sum q_i|$ vary in a family of dimension $g-5$, $N$
must vary in a family of dimension $2$, i.e., the points $t_1 $ and
$t_2$ are general in $C'$. Since $L$ is fixed this gives
\[
h^0 (\kappa^* g^1_2 -L )>0
\]
and $L$ has degree $4$.
\end{proof}

\subsection{}
If $d=4$, then $X$ is a Prym-embedded curve \cite{recillas74}. So $X$
deforms out of $\cJ_g$ into the locus of Prym varieties.

Let us then analyze the case $d=5$. Here $h^0 (L) =3$ so $g\leq 6$. By
the above, for $X$ to deform out of $\cJ_g$, it is necessary that, if
generically on a component of $Z(X)$ we have $h^0 (K -\sum q_i -L ) =
0$, then the image in $\cQ$ of the inverse image of that component in
$\tZ (X)$ does not generate $| I_2 (C) |$. Let $D_2\in X$ be such that
$h^0 (\sum q_i + D_2)\geq 2$. Since $L$ is in a $g^2_5$, any divisor
of $L$ spans a plane in $| K |^*$. Let us now distinguish the cases of
different genera.
\begin{description}
\item[g=4] Here $g-3 =1$ and $\sum q_i = q_1$. The variety $\cQ =
|I_2 (C) |$ is a point so for $X$ to deform out of $\cJ_4$ we need
that for all $q_1\in Z(X)$, $h^0 (K - L - q_1 ) > 0$. Any $g^2_5$ on
$C$ is of the form $| K- t |$ for some point $t$ on $C$. So $h^0 (K - L -
q_1 ) = h^0 (t - q_1 )$ is positive only when $t = q_1$. So for $X$ to
deform out of $\cJ_g$ we need $Z (X) = \{ t\}$. 

Let us now determine $Z(X)$. To say $h^0 ( D_2 + q_1 )\geq 2$ means of
course $| D_2 + q_1 |$ is one of the two possibly equal $g^1_3$ on
$C$. Denote this $g^1_3$ by $G$. So $D_2$ is also on
\[
X_2( G ) =\{ D_2: h^0 (G - D_2 ) > 0\}.
\]
First note that for any $g^1_3$ on $C$, $X_2 (g^1_3)\cong C$ is
irreducible and if $X$ contains it, then $Z (X) = C$ and $X$ cannot
deform. Next the intersection number of $X_2 (G)$ with $X$ is
\[
X\cdot X_2( G ) = (C_L -\delta )\cdot (C_{ G } -\delta ) = 15
- 5 - 3 - 3 = 4.
\]

We now find these four divisors of degree $2$ geometrically.  The
divisors of $g^2_5 = | K -t |\supset L$ are cut by planes through
$t$. A pencil of these planes whose base locus we denote by
$L_0\subset |K|^*, (L_0\cong\bP^1)$ cuts the divisors of $L$ on $C$
and the divisor $D_5$ of $L$ containing $D_2$ is cut by the span
$\langle L_0, t+D_2\rangle$ which is a plane. On the other hand, the
divisors of $G$ are cut on $C$ by a ruling $R$ of the unique quadric
$Q$ containing $C$. Since $X$ does not contain $X_2 (g^1_3)$ for any
$g^1_3$ on $C$, it follows that $L_0$ is not contained in $Q$. So
$L_0\cap Q$ is the union of two possibly equal points. Exactly one
line of $R$ passes through each of these points cutting two divisors
of $G$ on $C$. One of these divisors is the divisor of $R$ containing
$t$, say $E_2+t$ with $E_2\in X\cap X_2 (G)$. Writing the other
divisor as $t_1 +t_2 + t_3$, we have $t_i + t_j\in X\cap X_2 (G)$ for
all $i,j\in\{ 1,2,3\}$ which give us the other three points of $X\cap
X_2 (G)$. This means that $t_i\in Z(X)$ for all $i\in\{
1,2,3\}$. Therefore for $X$ to deform, we must have $t_1 = t_2 = t_3 =
t$. Therefore the two divisors of $R$ are equal to $3t$, in
particular, $L_0$ is tangent to $Q$.

If $C$ has another $g^1_3$ we repeat the above argument to obtain that
it is also equal to $|3t|$.

So we see that if $X$ deforms out of $\cJ_4$, then $C$ has only one
$g^1_3$ with a triple ramification point $t$ such that $5t\in L\subset
|K -t |$ and $X_2 (g^1_3)$ meets $X$ only at $2t$ with intersection
multiplicity $4$. Finally note that the facts $g^1_3 = |3t|$, $5t\in
L$ and $X$ does not contain $X_2 (g^1_3)$ imply that $L$ has no
base-points.

\item[g=5] Here $g-3 =2$ and $\sum q_i = q_1 + q_2$. The
linear system $|K - L| = | K -g^2_5 |$ is a $g^1_3$ on $C$, unique
because the genus is $\geq 5$. The variety $\cQ$ is a plane quintic
with a double point: it is the image of $C$ in $\bP^2 = | I_2 (C) |$
by the morphism associated to $g^2_5$. Every quadric of $\cQ$ has rank
$4$ except its double point $Q_0$ which has rank $3$. The singular
locus of $Q_0$ is a secant to $C$ and its ruling cuts the divisors of
$g^1_3$ on $C$. The intersection of the singular locus of $Q_0$ with
$C$ is the divisor $D_0$ such that $2g^1_3 + D_0\sim K$. The base
locus of $|I_2 (C) |$ in $|K|^*$ is the rational normal scroll traced
by the lines generated by the divisors of $g^1_3$.

To determine $Z(X)$ we first fix a general divisor $D_2\in C^{ (2) }$
and find all the divisors $q_1 + q_2 $ such that $h^0 ( D_2 + q_1 +
q_2)\geq 2$. To say $h^0 ( D_2 + q_ 1+ q_2 )\geq 2$ means $| D_2 + q_1
+ q_2|$ is a $g^1_4$ on $C$. To this $g^1_4$ is associated a quadric
of rank $4$ which then contains $\langle D_2\rangle$. Assuming $h^0
(g^1_3 - D_2 ) = 0$, there is exactly a pencil of quadrics of $| I_2
(C) |$ which contain $\langle D_2\rangle$. This pencil cuts $\cQ$ in
five points counted with multiplicities giving us five quadrics
counted with multiplicities, and for each quadric a choice of a ruling
containing $\langle D_2\rangle$. To each ruling is associated a
$g^1_4$ such that $h^0 ( g^1_4 -D_2 ) > 0$. These $g^1_4$ can be
described as follows. Assuming that $D_2\neq D_0$, there is a unique
divisor of $g^2_5$ which contains $D_2$. Let this divisor be $D_5$ and
write $D_5 = D_2 + s_1 + s_2 + s_3$. We have three $g^1_4$ containing
$D_2$ obtained as $|D_2 + s_i + s_j |$. Futhermore, if $D_2 = t_1 +
t_2$, we have two other $g^1_4$ containing $D_2$ obtained as $|g^1_3 +
t_i |$. It is not difficult to see that these are distinct for a
general choice of $D_2$.

Since $d\geq 4$, we can find $D_2\in X$ such that $h^0 (g^1_3 - D_2 )
= 0$. Taking such $D_2$ general in $X$ we can also assume $D_2\neq
D_0$. With the above notation, the possibly equal elements of $Z(X)$
that we obtain for $D_2$ are $s_i + s_j$ and $g^1_3 - t_i$. The last
two are contained in a divisor of $g^1_3 = |K-L|$, meaning they
satisfy $h^0 ( K-L -\sum q_i ) >0$. The pair $(D_2, s_i +s_j)\in\tZ
(X)$ is above $s_i+ s_j\in Z(X)$ and its image in $\cQ$ is the quadric
swept by the planes spanned by the divisors of $|D_2 +s_i +s_j
|$. This quadric is also the image of $s_k + g^1_3$ for $k\neq i,j$
since $s_k + g^1_3 = |K- D_2 -s_i -s_j|$. So the quadric is the image
of the point $s_k$ of $C$ in $\cQ$. Since the base divisor of $L$ has
degree at most $2$, for a general choice of $D_2$ as above, at least
one of the $s_i$ will be a general point of $C$, and as $D_2$ varies,
this point will trace all of $C$ and its image in $\cQ$ will trace all
of $\cQ$. So for $X$ to deform we also need $h^0 (K-L-s_i-s_j) = h^0
(g^1_3 -s_i-s_j) >0$ for all $i\neq j$. This implies $s_1+ s_2 +
s_3\in g^1_3$ and since the divisor $s_1 + s_2 + s_3$ is not fixed, we
obtain $L = g^1_3 + D_2$. This contradicts the generality of
$D_2$. Therefore $X$ cannot deform out of $\cJ_5$.

\item[g=6]
Here $g-3=3$ and $\sum q_i = q_1 + q_2 + q_3$. The curve $C$ is a
smooth plane quintic and $ K\sim 2 g^2_5\sim 2L$. The variety $Sing
(\Theta )$ is the image of $C\times C$ via $(p,q)\mapsto |g^2_5 -p+q
|$ (see e.g. \cite{ACGH} p. 264). So every complete $g^1_5$ on $C$ has
exactly one base point. Since $C$ embeds in $\bP^2$ by the map
associated to $g^2_5$, given $t_1 + t_2 = D_2\in X$, there is a unique
divisor $D_5 = D_2 + s_1 + s_2 + s_3$ of $g^2_5$ containing it. The
one-parameter family $Z(D_2)$ of $g^1_5$ such that $h^0 ( g^1_5 -D_2
)>0$ has six components: one component is the family of pencils in
$g^2_5$ passing through $D_5$, two components are families of complete
$g^1_5$ obtained as $|g^2_5 -t| +t_i$ where $t$ varies in $C$, and the
last three components are families of complete $g^1_5$ obtained as $|
D_2 + s_i + s_j + t |$ with $t$ varying in $C$. So altogether (and
counting multiplicities) $Z(D_2)$ is the union of a smooth rational
curve and $5$ copies of $C$. The divisors $\sum q_i$ for the rational
component are all equal to $s_1 + s_2 + s_3$ for which $h^0 ( K- L
-\sum q_i ) = h^0 ( g^2_5 - s_1 - s_2 - s_3) = h^0 ( D_2 ) > 0$. The
divisors $\sum q_i$ for the first two copies of $C$ in $Z(D_2)$ are
$g^2_5 - t - t_j$ so we see that they satisfy $h^0 ( K- L -\sum q_i )
= h^0 (g^2_5 - (g^2_5 - t - t_j )) = h^0 ( t + t_j )>0$. The divisors
$\sum q_i$ for the last three copies of $C$ in $Z(D_2)$ are $s_i +
s_j + t$ and so for $t$ general, we have $h^0 ( K- L -\sum q_i ) = h^0
(g^2_5 - s_i - s_j - t)= h^0 ( D_2 + s_k - t) = 0$. As we saw, here
$\sum q_i = s_i + s_j + t = g^2_5 - D_2 - s_k +t$. The divisors that
we obtain in $\tZ (X)$ map in $Sing (\T )$ to $g^2_5 - s_k + t$. As
$D_2$ varies in $X$, the points $s_k$ and $t$ vary freely in $C$ and
$g^2_5 - s_k + t$ traces all of $Sing (\T )$. So we see that $X$
cannot deform with $JC$ out of $\cJ_6$. Note however, that if we
degenerate the plane quintic to a singular one, then $X$ might deform.

\end{description}

\section{Appendix}\label{sectapp}

\subsection{The cohomology of some sheaves on $C^{ (n) }$}\label{LELE}

We calculate the cohomologies of some sheaves on $C^{ (n) }$
for an integer $n\geq 2$. Recall that $\pi_n : C^n\ra
C^{ (n) }$ is the natural morphism and let $\D_{i,j}^n$ ($1\leq i < j\leq
n$) be the diagonals of $C^n$. Also let $pr_i : C^n\ra C$ be the
$i$-th projection. Then
\[
\pi_n^* (\omega_{ C^{ (n) }})\cong\omega_{ C^n }\otimes\cO_{ C^n }
(\sum_{ 1\leq i< j\leq n } -\D_{ i,j }^n )
\]
by the Hurwitz formula, and
\[
\omega_{ C^n }\cong\otimes_{ i=1 }^n pr_i^*\omega_C.
\]
For any non-trivial divisor class $b$ of degree $0$ on $C$, the
intersection $\Theta .\Theta_b$ is easily seen to be reduced and its
inverse image in $C^{ (g-1) }$ is
\[
\{ D\in C^{ (g-1) } : h^0 ( D -b ) > 0\} =\{ D\in C^{ (g-1) } : h^0 (
K +b -D) > 0\}.
\]
If $n\leq g-1$, the restriction of this to $C^{ (n) }_{\sum_{ i=1 }^{
g-1-n } q_i}$ is
\[
\{ D\in C^{ (n) } : h^0 ( K +b -\sum_{ i=1 }^{ g-1-n } q_i - D) > 0\}
\]
whose pull-back to $C^n$ by $\pi_n$ is in the linear system
\[\begin{split}
|pr_1^*\cO_C (K + b -\sum_{ i=1 }^{ g-1-n } q_i )\otimes\ldots
\\
\ldots\otimes
pr_n^*\cO_C (K+ b -\sum_{ i=1 }^{ g-1-n } q_i)\otimes\cO_{ C^n }
(-\sum_{ 1\leq i < j\leq n }\D_{i,j}^n) |
\end{split}
\]
as can be easily seen by restricting to fibers of the various
projections $C^n\ra C^{ n-1}$ and using the See-saw Theorem.

More generally, for any divisor $E$ on $C$, let $\cL_{n,E}$ and $\cL_{
  n,E }'$ be the invertible sheaves on $C^{ (n) }$ whose inverse
images by $\pi_n$ are isomorphic to
\[
pr_1^*\cO_C (E )\otimes\ldots\otimes pr_n^*\cO_C (E)
\]
and
\[
pr_1^*\cO_C (E )\otimes\ldots\otimes pr_n^*\cO_C (E)\otimes\cO_{ C^n
  } (-\sum_{ 1\leq i < j\leq n }\D_{i,j}^n)
\]
respectively. We will calculate the cohomologies of $\cL_{n,E}$ and
$\cL_{ n, E}'$.

Since
\[
\pi_n^* \cL_{n,E}\cong pr_1^*\cO_C (E )\otimes\ldots\otimes
pr_n^*\cO_C (E),
\]
the sheaf $\cL_{n,E}$ is the invariant subsheaf of
\[
{\pi_n}_*\pi_n^* \cL_{n,E} = {\pi_n}_* (pr_1^*\cO_C (E
)\otimes\ldots\otimes pr_n^*\cO_C (E))
\]
for the action of the symmetric group $\Sg_n$. We claim that $\cL_{ n, 
  E}'$ is the skew-symmetric subsheaf of
\[
{\pi_n}_* (pr_1^*\cO_C (E )\otimes\ldots\otimes pr_n^*\cO_C (E))
\]
for the action of $\Sg_n$. To see this, first note that any
skew-symmetric local section of ${\pi_n}_* (pr_1^*\cO_C (E
)\otimes\ldots\otimes pr_n^*\cO_C (E))$ vanishes on
the diagonal $\D^n$ of $C^{ (n) }$. Conversely, pulling back to $C^n$,
we see that the ideal sheaf of any of the diagonals is generated by
skew-symmetric tensors.

Therefore the cohomology groups of $\cL_{ n,E }$ (resp. $\cL_{ n,E
}'$) are the invariant (resp. skew-symmetric) parts of the cohomology
groups of ${\pi_n}_* (pr_1^*\cO_C (E )\otimes\ldots\otimes pr_n^*\cO_C
(E ))$ under the action of $\Sg_n$. Or, equivalently, since $\pi_n$ is
finite, the invariant (resp. skew-symmetric) parts of the cohomology
groups of $pr_1^*\cO_C (E )\otimes\ldots\otimes pr_n^*\cO_C (E )$
under the action of $\Sg_n$.  By the K\"unneth formula the cohomology
groups of $pr_1^*\cO_C (E )\otimes\ldots\otimes pr_n^*\cO_C (E )$ are
\[
H^{i} (pr_1^*\cO_C (E )\otimes\ldots\otimes pr_n^*\cO_C (E
))\cong\sum H^{j_1} (E)\otimes\ldots\otimes H^{ j_n} (E) 
\]
where the $n$-uples $(j_1 ,\ldots , j_n )$ describe the set of
$n$-uples of elements of $\{0,1\}$, $i$ of which are equal to $1$ and
the rest equal to $0$. The action of $\Sg_n$ on each of these cohomology
groups is super-symmetric: for instance any transposition $\tau$
exchanging $l$ and $k$ sends an element $e_1\otimes\ldots\otimes
e_l\otimes\ldots\otimes e_k\otimes\ldots\otimes e_n$ to $(-1)^{j_l
  j_k} e_1\otimes \ldots\otimes e_k\otimes\ldots\otimes
e_l\otimes\ldots\otimes e_n$. From this one easily calculates that the
invariant parts of the cohomology groups are
\[
\begin{array}{l}
H^0 (\cL_{n,E })\cong S^n H^0 (E) \\
H^1 (\cL_{ n,E })\cong H^1 (E)\otimes S^{n-1 } H^0 (E)\\
H^2 (\cL_{ n,E })\cong\wedge^2 H^1 (E)\otimes S^{ n-2 } H^0 (E)\\
\hspace{20pt}\vdots\\
H^{ n-1} (\cL_{ n, E})\cong\wedge^{n-1 } H^1 (E)\otimes H^0 (E)\\
H^d (\cL_{ n, E} )\cong\wedge^n H^1 (E).
\end{array}
\]
Similarly, the skew-symmetric parts of the cohomology groups are
\[
\begin{array}{l}
H^0 (\cL_{n,E }')\cong\wedge^n H^0 (E) \\
H^1 (\cL_{ n,E }')\cong H^1 (E)\otimes\wedge^{n-1 } H^0 (E)\\
H^2 (\cL_{ n,E }')\cong S^2 H^1 (E)\otimes\wedge^{ n-2 } H^0 (E)\\
\hspace{20pt}\vdots\\
H^{ n-1} (\cL_{ n, E}')\cong S^{n-1 } H^1 (E)\otimes H^0 (E)\\
H^n (\cL_{ n, E}' )\cong S^n H^1 (E).
\end{array}
\]

\subsection{Useful exact sequences of cohomology groups}\label{ontoC2}

Let $a$ be such that $\Theta_a\supset C^{ (2) }$ and $Sing (\Theta_a
)\not\supset C^{ (2) }$. Then $-a +\sum p_i\sim\sum q_i$ and $h^0
(\sum q_i ) =1$. As we saw in \ref{LELE} we have $\cO_{ C^{
(2) }} (\Theta_a )\cong\cL_{ 2,K -\sum q_i}'$. Consider the composition
\[
C^{ (2) }\cong C_{\sum_{i=1}^{g-3} q_i }^{ (2) }\subset
C_{\sum_{i=1}^{ g-4 } q_i }^{ (3) }\subset\ldots\subset C_{ q_1 }^{
  (g-2) }\subset C^{ (g-1) }\ra\Theta_a\subset A
\]
For $3\leq n\leq g-1$, we have the exact sequence
\[\begin{split}
0\lra\cO_{ C^{ (n) }_{\sum_{i=1}^{g-1-n} q_i }} (\Theta_a - C^{ (n-1)
}_{\sum_{i=1}^{g-n} q_i })\lra\cO_{ C^{ (n) }_{\sum_{i=1}^{g-1-n} q_i
}} (\Theta_a )\lra\\
\lra\cO_{ C^{ (n-1) }_{\sum_{i=1}^{g-n} q_i }} (\Theta_a
)\lra 0.\end{split}
\]
For each $i$, by \ref{LELE}, we have
\[\begin{split}
H^i (\cO_{ C^{ (n) }_{\sum_{i=1}^{g-1-n} q_i }} (\Theta_a - C^{ (n-1)
}_{\sum_{i=1}^{g-n} q_i }))\cong \\
\cong S^i H^1 (K -\sum_{i=1}^{g-n} q_i
)\otimes\Lambda^{ n-i } H^0 (K -\sum_{i=1}^{g-n} q_i ),\end{split}
\]
\[
H^i (\cO_{ C^{ (n) }_{\sum_{i=1}^{g-1-n} q_i }} (\Theta_a ))\cong S^i
H^1 (K -\sum_{i=1}^{g-1-n} q_i )\otimes\Lambda^{ n-i } H^0 (K
-\sum_{i=1}^{g-1-n} q_i ),
\]
\[
H^i (\cO_{ C^{ (n-1) }_{\sum_{i=1}^{g-n} q_i }} (\Theta_a ))\cong S^i
H^1 (K -\sum_{i=1}^{g-n} q_i )\otimes\Lambda^{ n-1-i } H^0 (K
-\sum_{i=1}^{g-n} q_i )
\]
and the map on cohomology
\[
H^i (\cO_{ C^{ (n) }_{\sum_{i=1}^{g-1-n} q_i }} (\Theta_a - C^{ (n-1)
}_{\sum_{i=1}^{g-n} q_i }))\lra H^i (\cO_{ C^{ (n)
}_{\sum_{i=1}^{g-1-n} q_i }} (\Theta_a ))
\]
is obtained from the inclusion
\[
H^0 (K -\sum_{i=1}^{g-n} q_i )\inj H^0 (K -\sum_{i=1}^{g-1-n} q_i )
\]
(note that the dimension of $H^1 (K -\sum_{i=1}^{g-n} q_i )$ and $H^1
(K -\sum_{i=1}^{g-1-n} q_i )$ is $1$). It follows that for all $i$ the sequence
\[
0\lra H^i (\cO_{ C^{ (n) }_{\sum_{i=1}^{g-1-n} q_i }} (\Theta_a - C^{
(n-1) }_{\sum_{i=1}^{g-n} q_i }))\lra H^i (\cO_{ C^{ (n)
}_{\sum_{i=1}^{g-1-n} q_i }} (\Theta_a ))
\]
\[
\lra H^i (\cO_{ C^{ (n-1)
}_{\sum_{i=1}^{g-1-n} q_i }} (\Theta_a ))\lra 0
\]
is exact.
In particular, all the sections of $\cO_{ C^{ (2) } } (\Theta_a )$
vanish on $(Z_{ g-1 })_a\cap C^{ (2) }$, hence on $(Z_{ g-1 })_a\cap X$, so
they restrict to sections of $\cI_{ (Z_{ g-1 })_a\cap X} (\Theta_a )$ on
$X$.

\subsection{The cohomology of $\cL_{ 2,E} (-\Delta
)=\cO_{ C^{ (2) }} (C_E -\Delta ) = \cO_{ C^{ (2) }} (C_E -\Delta
)$}\label{appI2}

We use the exact sequence
\[
0\lra\cO_{ C^{ (2) }} (C_E  -\Delta )\lra \cO_{ C^{ (2) }} (C_E
)\lra\cO_{ C^{ (2) }} (C_E ) |_{\Delta }\lra 0\; .
\]
Under the isomorphism $\Delta\cong C$ we have $\cO_{ C^{ (2) }} (C_E )
|_{\Delta }\cong\cO_C (2E)$) and, by \ref{LELE}, we have the
long exact sequence of cohomology
\begin{equation}\label{2delta}\begin{split}
0\lra H^0 ( C^{ (2) }, C_E -2\delta )\lra S^2 H^0 (C, E)\lra H^0 (C,
2E)\lra
\\
\lra H^1 ( C^{ (2) }, C_E -2\delta )\lra H^0 (C, E)\otimes H^1 (C, E)\lra
\\
\lra H^1 (C , 2E)\lra H^2 ( C^{ (2) }, C_E
-2\delta )\lra\wedge^2 H^1 (C, E)\lra 0\; .\end{split}
\end{equation}
So, in particular, $H^0 (C^{ (2) }, C_E -2\delta )= I_2 (C, E)$ is the
space of quadratic forms vanishing on the image of $C$ in $|E|^*$ if $|E
|\neq\emptyset$.

Note that using our result in Appendix \ref{LELE} above, Pareschi and Popa
\cite{pareschipopa03} have computed the cohomology of $\cL_{ n, E }
(-\Delta )$ for $n > 2$ as well.


\providecommand{\bysame}{\leavevmode\hbox to3em{\hrulefill}\thinspace}

\end{document}